\documentclass[10pt]{article}\usepackage{amsmath,amssymb,amsthm}
\usepackage{amsfonts}
\usepackage{amsmath}

\usepackage[T1]{fontenc}
\usepackage[english]{babel}
\title{The randomization by Wishart laws and the Fisher information}
\author{G\'erard Letac\thanks{Institut de Math\'ematiques de Toulouse, 
Universit\'e Paul Sabatier, 31062 Toulouse, France. \texttt{gerard.letac@math.univ-toulouse.fr} }} 

\date{}\def \d{\frac{1}{2}\, }
\def\R{\mathbb{R}}
\def\P{\mathbb{P}}
\def\E{\mathbb{E}}
\def\tr{\,\textmd{tr}\,}\def\d{\frac{1}{2}\, }
\def\<{\langle}
\def\>{\rangle}

\begin{document}
\maketitle

\textit{To the memory of H\'el\`ene Massam, 1949-2020}

\begin{abstract} 
Consider  the centered  Gaussian  vector $X$  in $\R^n$ with  covariance matrix $ \Sigma.$ Randomize $\Sigma$ such that $
\Sigma^{-1}$ has a Wishart distribution with shape parameter $p>(n-1)/2$ and mean $p\sigma.$ We compute the density $f_{p,\sigma}$ of $X$ as well as the Fisher information $I_p(\sigma)$ of the model $(f_{p,\sigma} )$ when $\sigma $ is the parameter. For using the Cram\'er-Rao inequality, we also compute  the  inverse of $I_p(\sigma)$. The important point of this note is the fact that this inverse  is a linear combination of two simple operators on the space of symmetric matrices, namely 
$\P(\sigma)(s)=\sigma s \sigma$ and $(\sigma\otimes \sigma)(s)=\sigma \, \mathrm{trace}(\sigma s)$. The Fisher information itself is a linear combination $\P(\sigma^{-1})$ and $\sigma^{-1}\otimes \sigma^{-1}.$ Finally, by randomizing $\sigma $ itself, we make explicit the minoration of the second moments of an estimator of $\sigma$ by the Van Trees inequality:  here again, linear combinations of  $\P(u)$ and $u\otimes u$ appear in the results.

\end{abstract}
\vspace{4mm}\noindent{Keywords:} Moments of the Wishart distribution, Van Trees inequality, randomization of covariance. AMS classification: primary 62H10, secondary 62F15.

\section{Introduction}
 Let $n$ be  a fixed positive integer, and let $ E$ be  the Euclidean  space of symmetric real  matrices of order $n$. Let  $E_+\subset E$ be the cone of positive definite matrices.  Consider  the centered  Gaussian vector $X$  in $\R^n$ with  covariance matrix $ \Sigma\in E_+.$ Let us fix $\sigma\in E_+$ and $p>(n-1)/2.$ Suppose that $ \Sigma $ is itself random such that $U=\Sigma^{-1}$ has the Wishart distribution $\gamma_{p, \sigma}$, defined by its Laplace transform as follows: for all $s\in E_+$ we have
\begin{equation}\label{WISHLAPLACE}\int_{E_+}e^{-\tr(su)}\gamma_{p,\sigma}(du)=\frac{1}{\det (I_n+\sigma s)^p}.\end{equation} Since  $U\sim \gamma_{p,\sigma}$ and $X|U\sim N(0,U^{-1})$ the marginal density of $X$ is 
\begin{equation}\label{DENSITY}f_{p,\sigma}(x)=\frac{1}{(\sqrt{2\pi})^n}\int_{E_+}e^{-x^tux/2}(\det u)^{1/2}\gamma_{p,\sigma}(du).\end{equation}
We will give the explicit forms of $\gamma_{p,\sigma}(du)$ and of $f_{p,\sigma}(x)$ in a moment in  formulas  \eqref{WISHDENS} and \eqref{NORMDENS3}. We now consider $\sigma \mapsto f_{p,\sigma}(x)dx$ as a Fisher model and the aim of this note is to show that its Fisher information $I_p(\sigma)$ has an extremely simple form when using the good formalism. This good formalism uses 
two basic endomorphisms  $v\mapsto f(v)$ of the linear space $E$ . They depend on  a  parameter $u\in E$, namely $u\otimes u$ and $\P(u)$  and they are respectively defined  by \begin{eqnarray}\label{QUADRA1}
v&\mapsto& (u\otimes u)(v)=u\tr(uv)\\\label{QUADRA2}v&\mapsto& \P(u)(v)=uvu.\end{eqnarray}
The linear combinations of $u\otimes u$ and $\P(u)$ form a two dimensional linear space $L_u$ when $u\neq 0.$

Like $\sigma\otimes \sigma$ and $\P(\sigma)$, the Fisher information  $I_p(\sigma)$
is also a linear application of $E$ into itself . The aim of the present note is to prove  that there exist two numbers $a=a(p)$ and $b=b(p)$ such that
$$I_p(\sigma)=a\, \P(\sigma^{-1})+b\, \sigma^{-1}\otimes \sigma^{-1}$$ and to compute these numbers (formula \eqref{INFO2}). In other terms we are going to show that $I_p(\sigma)$ is in $L_{\sigma^{-1}}.$ The determinant of an element of $L_u$ being easy to compute (see Lemma 3.2 part 3), this leads to an explicit expression of the Jeffreys distribution (Corollary 3.3). 
In Lemma 3.2 part 1, we will see that an invertible element of $L_{\sigma^{-1}}$ is always in $L_{\sigma}.$ A consequence is that, in order to use the Cram\'er Rao inequality we have a simple expression of  $I_p(\sigma)^{-1},$ which is  in $L_{\sigma}$ (Proposition 4.1).

From  formula \eqref{INFO2}, it is now easy to randomize $\sigma$ by another  Wishart distribution and to apply the classical Van Trees inequality. Indeed, if the previously fixed $\sigma$ becomes the random  Wishart matrix $U$, with $ U\sim \gamma_{p_1, \sigma_1},$ then we use the fact   that $\E(U^{-1}\otimes U^{-1})$  and $\E(\P(U^{-1}))$ are well known to obtain  the explicit Van Trees minoration as an element of $L_{\sigma_1}.$

 Sections 2, 3 and 4 provide respectively $f_{p,\sigma},$ $I_p(\sigma)$ and its inverse. Section 5 gives general details about the Van Trees inequality, Section 6  presents the Van Trees minoration.  An appendix contains a proof of  the Van Trees inequality and  comments.

This note was the subject of an invited lecture at ISBA 2021 in a session in  the memory of H\'el\`ene Massam. I have been the coauthor of H\'el\`ene for several articles, including papers on the moments of the Wishart distributions. Several formulas obtained there are used here. 

\section
{The  randomized density $f_{p,\sigma}$ }

\vspace{4mm}\noindent \textbf{Proposition 2.1. }\begin{equation}\label{NORMDENS3}
f_{p,\sigma}(x)=\frac{1}{(2\pi)^{n/2}}\frac{\Gamma(p+\d)}{\Gamma(p-\frac{n-1}{2})}\frac{(\det \sigma)^{\d}}{(1+\d x^t\sigma x)^{p+\d}}.\end{equation}

\vspace{4mm}\noindent \textbf{Proof. } Recall first that the appropriate Gamma function for the real Wishart distribution is 
$$\Gamma_n(p)=(2\pi)^{n(n-1)/2}\prod_{j=1}^{n}\Gamma(p-\frac{j-1}{d}).$$ The linear space $E$ is considered as  a Euclidean space with the scalar product $\<u,v\>=\tr(uv)$ which is the trace of the product of the two symmetric matrices $u$ and $v.$ As a Euclidean space, its natural Lebesgue measure gives  mass one to any unit cube: here we part from a tradition (\cite{MUIRHEAD}) with would use a different normalization of the  Lebesgue measure on $E$. 
With these  notations the explicit form of the Wishart distribution $\gamma_{p,\sigma}$ for $p>(n-1)/2$, as defined by \eqref{WISHLAPLACE}, is 
\begin{equation}\label{WISHDENS}\gamma_{p,\sigma}(du)=e^{-\tr(\sigma^{-1}u)}
\frac{(\det u)^{p-\frac{n+1}{2}}}  {(\det \sigma)^p\Gamma_n(p)}
1_{E_+}(u)du,\end{equation} where $1_{E_+}$ is the indicator of the set $E_+$.
 As a consequence we express $f_{p,\sigma}(x)$ defined by \eqref{DENSITY} as follows
\begin{equation}\label{NORMDENS}
f_{p,\sigma}(x)=\frac{1}{(2\pi)^{n/2}}\int_{E_+}e^{-\<\d xx^t+\sigma^{-1},u\>}\frac{(\det u)^{p+\d-\frac{n+1}{2}}}{(\det \sigma)^p\Gamma_n(p)}du.\end{equation}
To compute \eqref{NORMDENS} we introduce the element $\sigma_1$ of $E_+$ defined by $\sigma_1^{-1}=\d xx^t+\sigma^{-1}$
Since the total mass of $\gamma_{p+\d,\sigma_1}(du)$ is one we can write
\begin{equation}\label{AUXI}\int_{E_+}e^{-\<\sigma_1^{-1},x\>}(\det u)^{p+\d-\frac{n+1}{2}}du=(\det\sigma_1)^{p+\d}\Gamma_n(p+\d).\end{equation}
Comparing \eqref{NORMDENS} and \eqref{AUXI} we get 
\begin{equation}\label{NORMDENS2}
f_{p,\sigma}(x)=\frac{1}{(2\pi)^{n/2}}\frac{(\det\sigma_1)^{p+\d}}{(\det \sigma)^p}\times \frac{\Gamma_n(p+\d)}{\Gamma_n(p)}\end{equation}
The last fraction is easy to compute from the definition of $\Gamma_n(p);$ we get \begin{equation}\label{RATIOGAMMA}\frac{\Gamma_n(p+\d)}{\Gamma_n(p)}=\frac{\Gamma(p+\d)}{\Gamma(p-\frac{n-1}{2})}.\end{equation} The calculation of $\frac{(\det\sigma_1)^{p+\d}}{(\det \sigma)^p}$ is harder. For this we introduce a simple linear algebra lemma:

\vspace{4mm}\noindent\textbf{Lemma 2.2.} Let $E$ be  a Euclidean space and $\mathrm{id}_E $ be the identity map on $E.$  Let also $a\in E$  and if $c$ is a real number. The following statements hold true:
\begin{itemize}
\item If $c\|a\|^2\neq 1$ then $(\mathrm{id}_E-c a\otimes a)^{-1}$ exists if and only if $c\|a\|^2\neq 1$  and is $\mathrm{id}_E+c\frac{a\otimes a}{1-c\|a\|^2}.$

\item $\mathrm{id}_E-ca\otimes a$ is positive definite if and only  $c\|a\|^2<1.$

\item $\det (\mathrm{id}_E-ca\otimes a)=1-c\|a\|^2.$ 
\end{itemize}

\vspace{4mm }\noindent\textbf{Proof.}  If $a\neq 0$  let $e_1=a/\|a\|$ and consider an orthonormal basis  $e=(e_1,\ldots,e_d)$ of $E$. Then the representative matrices  in the basis $e$ of $ a\otimes a$ and $\mathrm{id}_E-c a\otimes a$ are 

$$[a\otimes a]_e^e= \mathrm{diag}(\|a\|^2,0,\ldots,0),\ [\mathrm{id}_E-c a\otimes a]_e^e=\mathrm{diag}(1-c\|a\|^2,1,\ldots,1)$$
This leads immediately to the proof of the three claims.  \hfill
$ \square$

\vspace{4mm }\noindent We come back to the proof of Proposition 2.1.  From the definition of $\sigma_1$ we get
$$(\det \sigma_1)^{-1}=\det(\sigma^{-1}+\d xx^t)=(\det \sigma)^{-1/2}\det(I_n+\d  (\sigma^{1/2}x)(\sigma^{1/2}x)^t)(\det \sigma)^{-1/2}. $$
Apply Lemma 2.2 to $E=\R^n$, to $c=-1/2$ and to $a=\sigma^{1/2}x$
we get $$(\det \sigma_1)^{-1}=\det (\sigma^{-1})(1+\d \|\sigma^{1/2}x\|^2)=\det (\sigma^{-1})(1+\d x^t\sigma x).$$
Gathering this with \eqref{RATIOGAMMA} the expression \eqref{NORMDENS2} gives \eqref{NORMDENS3}. \hfill
$ \square$

\vspace{4mm}\noindent\textbf{Corollary 2.3.}  The \textit{a posteriori} distribution of $\Sigma^{-1}|X$ is $\gamma_{\sigma_1, p+\d}$ where $\sigma_1^{-1}=\sigma^{-1}+\d xx^t.$

\vspace{4mm}\noindent\textbf{Proof.} This is clear from the joint distribution of $(X,\Sigma^{-1})$ and the explicit form \eqref{NORMDENS3} of $f_{p,\sigma}.$\hfill
$ \square$

\section{The Fisher information $I_p(\sigma)$}
The next step is the calculation of the Fisher information  $I_p(\sigma)$ of the Fisher model $$\{f_{p,\sigma}(x)dx; \sigma \in E_+\}$$ on $\R^n$. It is better  not to call it the Fisher information matrix, since its representation by a matrix as an endomorphism of  the  space $E$ of symmetric matrices  is not appropriate.  The main result of this note is the following proposition, showing that $I_p(\sigma)$ is in the linear space $L_{\sigma^{-1}}:$

\vspace{4mm}\noindent \textbf{Proposition 3.1.}\begin{equation}\label{INFO2}I_p(\sigma)=\frac{1}{2(2p+3)}\left((2p+1)\P(\sigma^{-1})-\sigma^{-1}\otimes \sigma^{-1}\right)\end{equation}

\vspace{4mm}\noindent \textbf{Proof.}
From \eqref{NORMDENS3} we can write $$f_{p,\sigma}(x)=\frac{1}{(2\pi)^{n/2}}\frac{\Gamma(p+\d)}{\Gamma(p-\frac{n-1}{2})}
e^{\ell_x(\sigma)},$$ where $$\ell_x(\sigma)=\d\log \det \sigma-\left(p+\d\right)\log \left(1+\d x^t\sigma x\right).$$
Consider now for $x$ fixed in $\R^n$ the differential (or gradient) $\ell'_x(\sigma)$ of the real function $\sigma \mapsto \ell_x(\sigma)$ defined on $E_+\subset E$
and its second differential $\ell''_x(\sigma).$ We look for $$I_p(\sigma)=\int_{\R^n}\ell'_x(\sigma )\otimes \ell'_x(\sigma )f_{p,\sigma}(x)dx=-\int_{\R^n}\ell''_x(\sigma ) f_{p,\sigma}(x)dx$$
(recall that the notations $u\otimes u$ and $\P(u)$ for $u\in E$ are defined in \eqref{QUADRA1} and \eqref{QUADRA2}). To compute
$\ell'_x(\sigma )$ and $ \ell''_x(\sigma )$ we use the following three facts (\cite{LETACMASSAM1998}): \begin{itemize}\item the differential of $\sigma \mapsto \log \det \sigma $ is $\sigma^{-1}$; \item the differential of $\sigma \mapsto  \sigma^{-1} $ is $-\P(\sigma^{-1})$; 
\item the differential of $\sigma \mapsto x^t\sigma x=\<\sigma,xx^t\>$ is $xx^t.$

\end{itemize}
As a consequence we get
\begin{eqnarray*}
\ell'_x(\sigma )&=&\d \sigma^{-1}-\left(p+\d\right)\d\frac{xx^t}{1+\d x^t\sigma x}\\
\ell''_x(\sigma )&=&-\d \P(\sigma^{-1})+\left(p+\d\right)\frac{1}{4}\frac{xx^t\otimes xx^t}{\left(1+\d x^t\sigma x\right)^2}.
\end{eqnarray*}
We now introduce the endomorphism $J_p(\sigma)$ of $E$ defined by the integral 
\begin{eqnarray}\nonumber
J_p(\sigma)&=&\int_{\R^n}\frac{xx^t\otimes xx^t}{\left(1+\d x^t\sigma x\right)^2}f_{p,\sigma}(x)dx\\\nonumber&=&\frac{1}{(2\pi)^{n/2}}\frac{\Gamma(p+\d)}{\Gamma(p-\frac{n-1}{2})}\int_{\R^n}\frac{xx^t\otimes xx^t}{\left(1+\d x^t\sigma x\right)^2}\frac{(\det \sigma)^{\d}}{(1+\d x^t\sigma x)^{p+\d}}dx\\&=&\frac{(\det \sigma)^{\d}}{(2\pi)^{n/2}}\frac{\Gamma(p+\d)}{\Gamma(p-\frac{n-1}{2})}\int_{\R^n}\frac{xx^t\otimes xx^t}{(1+\d x^t\sigma x)^{p+2+\d}}dx\label{INTJ}
\end{eqnarray}
Its calculation is a crucial point since \begin{equation}\label{INFO}I_p(\sigma)=\d \P(\sigma^{-1})-\left(p+\d\right)\frac{1}{4}J_p(\sigma).\end{equation} We use for $a,b>0$ the integral representation 
$$\frac{1}{a^b}=\int_{0}^{\infty}e^{-ra}r^{b-1}\frac{dr}{\Gamma(b)}$$ applied to $a=1+\d x^t\sigma x$ and $b=p+2+\d.$ We get 

\begin{eqnarray*}J_p(\sigma)&=&\frac{\Gamma(p+\d)}{\Gamma(p-\frac{n-1}{2})}\times \frac{(\det \sigma)^{\d}}{(2\pi)^{n/2}}\int_{\R^n}(xx^t\otimes xx^t)\left(\int_0^{\infty}e^{-r(1+\d x^t\sigma x)}r^{p+\frac{3}{2}}\frac{dr}{\Gamma(p+\frac{5}{2})}\right)dx\\
&=&\frac{\Gamma(p+\d)}{\Gamma(p-\frac{n-1}{2})}\times\int_0^{\infty}e^{-r}r^{p+\frac{3}{2}}K(r,\sigma)\frac{dr}{\Gamma(p+\frac{5}{2})}
\end{eqnarray*}
where $K(r,\sigma)$ is the symmetric endomorphism of $E$ defined by the integral
$$K(r,\sigma)=\frac{(\det \sigma)^{\d}}{(2\pi)^{n/2}}\int_{\R^n}(xx^t\otimes xx^t)e^{-\d rx^t\sigma x}dx.$$
If $s\in E_+$ observe that 
$$\frac{(\det \sigma)^{\d}}{(2\pi)^{n/2}}\int_{\R^n}e^{-\<s,\d xx^t\>}e^{-\d rx^t\sigma x}dx=\frac{(\det \sigma)^{1/2}}{\det(s+r \sigma)^{1/2}}=\frac{1}{r^{n/2}\det(I_n+\frac{\sigma^{-1}}{r}s)^{1/2}}.$$
In other terms, using the definition of the Wishart distribution by its Laplace transform \eqref{WISHLAPLACE}, if $$U\sim \gamma_{\d,\frac{1}{r}\sigma^{-1}}$$
one can say that  if $X\sim N(0,\frac{1}{r}\sigma^{-1})$ then $U=XX^t/2\sim \gamma_{\d,\frac{1}{r}\sigma^{-1}}.$ This implies that
$$K(r,\sigma)=\frac{4}{r^{n/2}}\E(U\otimes U).$$ We are now in position to use Proposition 3.2 of \cite{LETACMASSAM1998} which says that if $U\sim \gamma_{p,\sigma}$ then
\begin{eqnarray}\label{HELENE1}\E(U\otimes U)&=&p^2\sigma\otimes \sigma+p\, \P(\sigma)\\
\label{HELENE2}\E(\P(U))&=&\d p\, \sigma\otimes \sigma+(\d p+p^2)\P(\sigma)\end{eqnarray}
Applying \eqref{HELENE1} to the pair $ (\d,\frac{1}{r}\sigma^{-1})$ we get 

$$\E(U\otimes U)=\frac{1}{4r^2}\sigma^{-1}\otimes \sigma^{-1}+\frac{1}{2r^2} \P(\sigma^{-1}),$$
leading to the following expression for  $K(r,\sigma)$ and  $J_p(\sigma).$

\begin{eqnarray*}K(r,\sigma)&=&\frac{4}{r^{\d n+2}}\left(\frac{1}{4}\sigma^{-1}\otimes \sigma^{-1}+\d \P(\sigma^{-1})\right),\\
J_p(\sigma)&=&4\left(\frac{\Gamma(p+\d)}{\Gamma(p-\frac{n-1}{2})\Gamma(p+\frac{5}{2})}\int_0^{\infty}e^{-r}r^{p+\frac{1-n}{2}-1}dr\right)\left(\frac{1}{4}\sigma^{-1}\otimes \sigma^{-1}+\d \P(\sigma^{-1})\right)\\&=&4\left(\frac{\Gamma(p+\d)\Gamma(p-\frac{n-1}{2})}{\Gamma(p-\frac{n-1}{2})\Gamma(p+\frac{5}{2})}\right)\left(\frac{1}{4}\sigma^{-1}\otimes \sigma^{-1}+\d \P(\sigma^{-1})\right)\\&=&4\frac{1}{(p+\frac{3}{2})(p+\d)}\left(\frac{1}{4}\sigma^{-1}\otimes \sigma^{-1}+\d \P(\sigma^{-1})\right).
\end{eqnarray*}
We can now  compute $I_p(\sigma).$ From the formula \eqref{INFO} we get \eqref{INFO2}, which concludes the argument. \hfill
$ \square$

\vspace{4mm}\noindent Before stating a corollary of Proposition 3.1, we need a lemma:

\vspace{4mm}\noindent\textbf{Lemma 3.2.} Let $E$ be the Euclidean space of symmetric matrices of order $n$ and let $E_+\subset E$ be the set of the positive definite ones. Fix $\sigma$ in $E_+$ and $c$ a real number. Then
\begin{itemize} \item $(\P(\sigma)-c\sigma\otimes \sigma)^{-1}$ exists if and only if $c\neq 1/n$ and is $\P(\sigma)+\frac{c}{1-nc}\sigma\otimes \sigma.$

\item $\P(\sigma)-c\sigma\otimes \sigma$ is positive definite if and only if $c<1/n.$
\item $\det(\P(\sigma)-c\sigma\otimes \sigma)=(\det\sigma)^{n+1}(1-cn).$

\end{itemize}
\vspace{4mm}\noindent\textbf{Proof.} Let us observe that $\P(\sigma^{-1/2})(\sigma)=I_n$ and that $\|I_n\|^2=\tr(I_n)=n.$ Therefore, 
\begin{eqnarray*}\P(\sigma)-c\sigma\otimes \sigma&=& \P(\sigma^{1/2})\left(\mathrm{id}_E-c(\P(\sigma^{1/2})\sigma)\otimes (\P(\sigma^{1/2})\sigma)\right)\P(\sigma^{1/2})\\&=&\P(\sigma^{1/2})\left(\mathrm{id}_E-cI_n\otimes I_n\right)\P(\sigma^{1/2})\end{eqnarray*}
Applying Lemma 2.2  to $a=I_n$ shows the first two points as well as $\det \left(\mathrm{id}_E-cI_n\otimes I_n\right)=1-cn.$ For proving that $\det \P(\sigma)=(\det\sigma)^{n+1}$ let us assume without loss of generality that $\sigma=\mathrm{diag}(\lambda_1,\ldots,\lambda_n).$ This implies that if $v\in E$ we have \begin{equation}\label{EXPLICIT}\P(\sigma)(v)=\sigma v\sigma =(\lambda_i\lambda_jv_{ij})_{1\leq i,j\leq n}.\end{equation} If $(e_1,\ldots,e_n)$ is the canonical basis of $\R^n,$ we create a natural basis $f$ for $E$ made with the matrices $e_{ii}=e_i\otimes e_i$ for $i=1,\ldots,n$ and $e_{ij}=e_i\otimes e_j+ e_j\otimes e_i$ for $1\leq i<j\leq n.$  An interpretation of \eqref{EXPLICIT} is that $\P(\sigma)(e_{ij})=\lambda_i\lambda_je_{ij}$ for $1\leq i\leq j\leq n$, or that $e_{ij}$ is an eigenvector of $\P(\sigma)$ for the eigenvalue $\lambda_i\lambda_j.$ Since $(e_{ij})_{1\leq i\leq j\leq n}$ is a basis of $E$ we can claim that 
$$\det \P(\sigma)=\prod_{1\leq i\leq j\leq n}\lambda_i\lambda_j
=(\lambda_1\ldots\lambda_n)^{n+1}=(\det \sigma)^{n+1},$$ which concludes the proof.  \hfill
$ \square$

\vspace{4mm}\noindent\textbf{Corollary 3.3.} For the Fisher model $\{f_{p,\sigma}(x)dx; \sigma \in E_+\}$, the Jeffreys   \textit{a priori} measure  $\frac{1}{\sqrt{\det I_p(\sigma)}}1_{E_+}(\sigma)d\sigma$  is given by
$$\det I_p(\sigma)=\left(\frac{2p+1}{2(2p+3)}\right)^{n(n+1)/2}\left(1-\frac{n}{2p+1}\right)\, (\det \sigma)^{-n-1}.$$

\vspace{4mm}\noindent\textbf{Proof.} We need only  to apply Lemma 3.2 to $c=1/(2p+1)$  and  replace $\sigma$ by $\sigma^{-1}.$ Recall that $\dim E=n(n+1)/2.$ Note that $p>(n-1)/2$, or $ 1/(2p+1)<1/n.$ \hfill
$ \square$

\vspace{4mm}\noindent\textbf{Remarks.}
\begin{itemize} \item From the very definition of $I_p(\sigma)$ this symmetric endomorphism of the Euclidean space $E$ must be positive definite. Lemma 3.2 gives a direct proof of the fact that $I_p(\sigma)$ given by formula \eqref{INFO2} is positive definite.
\item To insist on the ubiquity of the symbols $\P(u)$  and $u\otimes u$ consider the simple problem of the calculation of the Fisher information $F_1(u)$ and $F_2(v)$ for  the two  Fisher models
$\{N(0,u^{-1})(dx), u\in E_+\}$  and $\{N(0,v)(dx), v\in E_+\}$ . These models are  unconcerned with the randomization of $u$ by Wishart law $\gamma_{p,\sigma}.$ In the first case the log likelihood is $\ell_x(u)=\d x^tux+\d \log \det u.$ This leads  to 
$$\ell'_x(u)=\d xx^t+\d u^{-1},\ \ell'_x(u)=-\d\P(u^{-1}), \ F_1(u)=\d\P(u^{-1}).$$ Instead of a  direct calculation, it is shorter to deduce $F_2(v)$ from $F_1(u)$ by the change of variable formula for Fisher information by $u=v^{-1}.$ Since the differential of $v\mapsto v^{-1}$ is $-\P(v^{-1})$ we get 
$$F_2(v)=(-\P(v^{-1}))^T\circ F_1(v^{-1})\circ (-\P(v^{-1}))=\d \P(v_{-1})\P(v)\P(v^{-1})=\d\P(v^{-1}).$$ Of course this direct proof can be avoided by observing that the model $\{N(0,u^{-1})(dx), u\in E_+\}$ is an exponential family.
\item It can be noted that in formula \eqref{INFO2} the dimension $n$  does not appear in the coefficients. 
\end{itemize}
\section{Explicit form of the Cram\'er- Rao inequality}
Consider a  random variable $\Sigma$   valued in $E$ as  an unbiased estimator of $\sigma$ for the Fisher model $\{f_{p,\sigma}(x)dx; \sigma\in E_+\}$ on $\R^n.$  The Cram\'er- Rao inequality says that $$\E((\Sigma-\sigma)\otimes (\Sigma-\sigma))- I_p(\sigma)^{-1}$$ is semi positive definite, where $I_p(\sigma)$ has been computed in \eqref{INFO2}. We show that $I_p(\sigma)^{-1}$ is in the space $L_{\sigma}$  and is quite explicit:

\vspace{4mm}\noindent\textbf{Proposition 4.1.} 

$$I_p(\sigma)^{-1}=\frac{2(2p+3)}{2p+1}\left(\P(\sigma)+\frac{1}{2p+1-n}\sigma\otimes \sigma\right).$$

\vspace{4mm}\noindent \textbf{Proof.} It follows from the formula \eqref{INFO2} and from Lemma 3.2 where $\sigma$ is replaced by $\sigma^{-1}$ and where $c=1/(2p+1)<1/n.$  \hfill
$ \square$

\section{The density information and the Van Trees inequality}
In this section we describe the Van Trees inequality discovered in \cite{vanTREES}. To this end, we define first the density information. Consider the probability density $\theta\mapsto \lambda(\theta)$ defined on an open set $\Theta\subset E$ where $E$ is a Euclidean space. We assume that

\begin{itemize}\item $\lambda$ is strictly positive on $\Theta$;\item the differential $\theta\mapsto \lambda'(\theta)$ exists and is continuous; \item $\lambda$ is zero at the boundary of $\Theta$  and vanishes at infinity. More specifically, for any $\epsilon >0$ there exists a compact $K_{\epsilon}\subset \Theta$ such that $\lambda (\theta)\leq \epsilon$ for all $\theta\in \Theta\setminus K_{\epsilon}.$ \end{itemize}
As a consequence $g=\log \lambda$ is well defined. The density information $I_{\lambda}$ is the symmetric endomorphism of $E$ defined by 
$$I_{\lambda}=\int_{\Omega}(g'(\theta)\otimes g'(\theta))\lambda(\theta)d\theta.$$ In the particular case $\Theta=E$ this density information could be considered as the Fisher information of the location model $\{\lambda(\theta-t)d\theta; t\in E\}.$ However the result

\begin{equation}\label{DENSINFO2}I_{\lambda}=-\int_{\Omega}g''(\theta)\lambda(\theta)d\theta.\end{equation}
is not proved as easily as its analog for the Fisher information. A necessary  condition for applying the Stokes formula is that $g'(\theta)\lambda (\theta)$ is zero on the boundary $\partial \Theta$ so that we can write $0=\int_{\Theta}(g'(\theta) \lambda(\theta))'d\theta.$

Independently, consider the Fisher model $P= \{P_{\theta}(dw); \theta\in \Theta\}$ defined on some measured space $(\Omega,\nu)$. For convenience we write $$P_{\theta}(dw)=e^{\ell_w(\theta)}\nu(dw)$$ and we consider its Fisher information $I_P(\theta)$ which is the symmetric endomorphism of $E$ equal to 
$$I_P(\theta)=\int_{\Omega}\ell'_w(\theta)\otimes \ell'_w(\theta)P_\theta(dw)=-\int_{\Omega}\ell''_w(\theta)P_\theta(dw).$$    

The Van Trees inequality is a Bayesian result. Let $X(w)$ be any estimator of $\theta$ -not necessarily unbiased. Suppose now that $\theta$ is randomized by $\lambda(\theta)$.
 What is the generalization of the Cram\'er-Rao theorem to this situation? That is: what is a minorant of the endomorphism of $E$ equal to $$C=\E((X-\theta)\otimes (X-\theta))$$ (in the sense of the partial ordering $A\geq B$ of the symmetric endomorphisms of $E$ defined by $B-A$ semi positive definite)? Here, the expectation  $C$ is taken on $\Omega\times \Theta$ with respect to the probability 
$P_{\theta}(dw)\lambda(\theta)d\theta$. The Van Trees inequality says that if $$D=I_{\lambda}+\int_{\Theta}I_P(\theta)\lambda(\theta)d\theta,$$
then $C\geq D^{-1}$, \textit{i.e.} $C- D^{-1}$ is semi positive definite. We give a proof in the appendix, with some comments explaining the differences with an ordinary Bayesian version  of the Cram\'er- Rao inequality. 

\section{Randomizing $I_p(\sigma)$ by  a Wishart distribution}

In this section, we apply the Van Trees inequality to the Fisher model $$\{P_{u}(dx); u\in E_+\}=\{f_{p,u}(x)dx; u\in E_+\}$$ (for clarity in the sequel, we have  replaced $\sigma$ or $\theta$ by $ u\in E_+$) and the density $\lambda(u)$ of $U\sim \gamma_{p_1,\sigma_1}$ and we will show that the Van Trees minoration is in $L_{\sigma_1}.$ 
 To perform this program we shall need the following formulas when $U\sim \gamma_{p,\sigma}:$ with $p>\frac{n+1}{2}$ for \eqref{WISHINV} and 
$p>\frac{n+3}{2}$ for \eqref{WISHINV2} and  \eqref{WISHINV3}:

\begin{eqnarray}
\label{WISHINV}\E(U^{-1})&=&\frac{\sigma^{-1}}{p-\frac{n+1}{2}};\\
\label{WISHINV2}\E(U^{-1}\otimes U^{-1})&=& \frac{1}{(p-\frac{n+3}{2})(p-\frac{n+1}{2})(p-\frac{n}{2})}\left((p-\frac{n+2}{2})\sigma^{-1}\otimes \sigma^{-1}+\P(\sigma^{-1})\right);\\
\label{WISHINV3}\E(\P(U^{-1}))&=& \frac{1}{(p-\frac{n+3}{2})(p-\frac{n+1}{2})(p-\frac{n}{2})}\left(\d\sigma^{-1}\otimes \sigma^{-1}+(p-\frac{n+1}{2})\P(\sigma^{-1})\right).
\end{eqnarray}
Formula \eqref{WISHINV} is classical and can be found in \cite {MUIRHEAD} page 97 or \cite{LETACMASSAM1998} formula 6.1. Formula \eqref{WISHINV2} is  due to von Rosen \cite{vonROSEN} and  formula \eqref{WISHINV2} is due to Das Gupta \cite{DASGUPTA}. A good synthesis is \cite{MATSUMOTO}.  Generalisations to symmetric cones appear in \cite{LETACMASSAM2000} page 137 and \cite{LETACMASSAM2004} page 307. 

  \vspace{4mm}\noindent\textbf{Proposition 6.1.}
Let $U\sim \gamma_{p_1,\sigma_1}$ with $p_1-\frac{n+3}{2}>0.$ and density $\lambda(u).$  Let $I_p(\sigma)$ as in \eqref{INFO2}.   
Then there exists numbers $a=a(p,p_1)$ and $b=b(p,p_1)$ such that
$$\left(I_{\lambda}+\int_{E_+}I_p(u)\lambda(u)du)\right)^{-1}=a\P(\sigma_1)+b\sigma_1\otimes \sigma_1.$$

\vspace{4mm}\noindent\textbf{Proof.} We begin with the calculation of the density information $I_\lambda.$ We consider \begin{eqnarray}\nonumber g(u)&=&\log \lambda(u) =-\<u,\sigma_1^{-1}\>+(p_1-\frac{n+1}{2})\log( \det u)+ \mathrm{constant}\\\label{ICASEA} g'(u)&=&-\sigma_1^{-1}+(p_1-\frac{n+1}{2})u^{-1},\ \ g''(u)=-(p_1-\frac{n+1}{2})\P(u^{-1}).\end{eqnarray}
First we have to insure that we can use \eqref{DENSINFO2}, by checking that if $u\to u_0\in \partial E_+$ then $g'(u)\lambda(u)$ tends to 0. If $\lambda_1,\ldots,\lambda_n$ are the eigenvalues of $u$ and $\lambda^0_1,\ldots,\lambda^0_n$ are the eigenvalues of $u_0$ expressed in a decreasing order, we can claim that $\lambda_i\to _{u\to u_0}\lambda_i^0.$ In particular $m=\min_{i}\lambda_i\to 0$  and there exists $M>1$ such that $\lambda_i<M$ for all $i=1,\ldots,n$ and therefore $\det u<mM^{n}.$  Furthermore the Frobenius norm of $u^{-1}$ satisfies  $$\|u^{-1}\|^2=\frac{1}{\lambda_1^2}+\cdots+\frac{1}{\lambda_n^2}\leq \frac{n}{m^2}.
$$The fact that $p_1>\frac{n+3}{2}$ implies that $p_1-\frac{n+1}{2}=1+\epsilon>1$. In particular $\lambda(u)\to _{u\to u_0} 0$ and 
$$\|u^{-1}\|(\det u)^{p_1-\frac{n+1}{2}}\leq \frac{\sqrt{n}}{m}m^{1+\epsilon}M^{n(1+\epsilon)}\to_{u\to u_0} 0.$$
Finally 
$$\|g'(u)\lambda(u)\|\leq  \|\sigma_1^{-1}\|\lambda(u)+(p_1-\frac{n+1}{2})\|u^{-1}\|\lambda(u)\to_{u\to u_0} 0.$$

In the sequel, to avoid unnecessarily  complicated formulas we use the  notation
$$A_i=p_1-\frac{n+i}{2}.$$ We can use \eqref{DENSINFO2} here and use \eqref{WISHINV3} for computing $I_{\lambda}$ from \eqref{ICASEA}.  We get 
\begin{equation}\label{ICASEA2}I_{\lambda}=\frac{1}{A_3A_0}\left(A_1\P(\sigma_1^{-1})+\d\sigma_1^{-1}\otimes \sigma_1^{-1}\right).\end{equation}

The second ingredient  of the Van Trees inequality is $\int_{E_+}I_p(u)\lambda(u)du.$ Using \eqref{INFO2}:
\begin{eqnarray}\nonumber
\int_{E_+}I_p(u)\lambda(u)du&=& \frac{2p+1}{2(2p+3)}\int_{E_+}\P(u^{-1})\gamma_{p_1,\sigma_1}(du)-\frac{1}{2(2p+3)}\int_{E_+}(u^{-1}\otimes u^{-1})\gamma_{p_1,\sigma_1}(du)\\\nonumber&=&\frac{2p+1}{2(2p+3)}\times \frac{1}{A_3A_0}P(\sigma_1^{-1})+\frac{2p+1}{2(2p+3)}\times \frac{1}{2A_3A_1A_0}(\sigma_1^{-1}\otimes \sigma_1^{-1})\\\nonumber
&&-\frac{1}{2(2p+3)}\times \frac{1}{A_3A_1A_0}P(\sigma_1^{-1})-\frac{1}{2(2p+3)}\times \frac{A_2}{A_3A_1A_0}(\sigma_1^{-1}\otimes \sigma_1^{-1})\\\label{ICASEA3}&=&\frac{1}{2(2p+3)A_3A_1A_0}\left((2pA_1+A_3)\P(\sigma_1^{-1})-(A_3-p)(\sigma_1^{-1}\otimes \sigma_1^{-1})\right).
\end{eqnarray}

Finally we get $I_{\lambda}+\int_{E_+}I_p(u)\lambda(u)du=A\P(\sigma_1^{-1})-B(\sigma_1^{-1}\otimes \sigma_1^{-1})$
where the real coefficients $A$ and $B$ are obtained from \eqref{ICASEA2} and \eqref{ICASEA3}. Applying now Lemma 4.2 we get finally 
the minoration of the Van Trees inequality, namely  
$$\left(I_{\lambda}+\int_{E_+}I_p(u)\lambda(u)du)\right)^{-1}=\frac{\P(\sigma_1)}{A}+\frac{B}{A(A+nB)}(\sigma_1\otimes \sigma_1).$$ \hfill
$ \square$

\section{Appendix: proof of the Van Trees inequality}  \vspace{4mm}\noindent \textbf{Proposition 7.1.} Let $(P_{\theta})_{\theta\in \Theta}$ be a regular Fisher model on $\Omega$  where  $\Theta$ is an open subset of the Euclidean space $E=\R^k$.  Denote its Fisher information by $I(\theta).$ Let $\lambda(\theta)$ be  a continuous  probability density on   $\overline{\Theta}$ which is continuously differentiable and positive  on $\Theta$,  and which is zero on the boundary $\partial \Theta=\overline{\Theta}\setminus \Theta$ and  at infinity: that means  that for any $\epsilon>0$ there exists a compact subset $K_{\epsilon}$ of $\Theta$ such that $\lambda(\theta)\leq \epsilon$ for $\theta\in \Theta\setminus K_{\epsilon}$.  Denote its density  information by $I_{\lambda}.$ Let $\Omega \ni w\mapsto X(w)\in E $  be an arbitrary estimator of $\theta.$ Under these circumstances the following $2k$ symmetric matrix 
\begin{equation}\label{vT}\left[\begin{array}{cc}\int_{\Omega\times \Theta}[(X(w)-\theta)\otimes (X(w)-\theta)]P_{\theta}(dw)\lambda(\theta)d\theta&I_k\\I_k&I_{\lambda}+\int_{ \Theta}I(\theta)\lambda(\theta)d\theta\end{array}\right]\end{equation} is semi positive definite. In particular the following fact, called the \textit{Van Trees inequality }holds: the $k$ symmetric matrix 
\begin{equation}\label{vT'}\int_{\Omega\times \Theta}[(X(w)-\theta)\otimes (X(w)-\theta)]P_{\theta}(dw)\lambda(\theta)d\theta-[I_{\lambda}+\int_{ \Theta}I(\theta)\lambda(\theta)d\theta]^{-1}\end{equation}
is semi  positive definite provided that the matrix  $I_{\lambda}+\int_{ \Theta}I(\theta)\lambda(\theta)d\theta$ is invertible. 

\vspace{4mm}\noindent \textbf{Proof.} Denote $P_{\theta}(dw)=e^{\ell_w(\theta)}\nu(dw)$ as usual and  write $\lambda(\theta)=e^{g(\theta)}.$ The basic tool of the proof is  Stokes' theorem, which says that if $f$ is a  function on $\overline{\Theta}$ valued in $\mathbb{R}^m$ which is sufficiently regular and which is zero on $\partial \Theta$ then $$\int_{\Theta}f'(\theta)d\theta=0.$$ We apply this principle to the two following functions 

\begin{enumerate}
\item  $f(\theta)=\lambda(\theta)e^{\ell_w(\theta)}$ and thus $m=1.$ Here \begin{equation}\label{M}f'(\theta)=     (g'(\theta)+\ell'_w(\theta))\lambda(\theta)e^{\ell_w(\theta)}  \end{equation}

\item $f(\theta)=\theta \lambda(\theta)e^{\ell_w(\theta)}$ and thus $m=\dim E.$ Here 
\begin{equation}\label{N}f'(\theta)=  [I_k+\theta\otimes(g'(\theta)+\ell'_w(\theta))]\lambda(\theta)e^{\ell_w(\theta)} \end{equation} 
\end{enumerate}
From (\ref{M}) and Stokes we get that $X(w)\otimes \int _{\Theta}(g'(\theta)+\ell'_w(\theta))\lambda(\theta)e^{\ell_w(\theta)}d\theta=0$ and thus integrating with respect to $\nu(dw)$ we have 
\begin{equation}\label{M'}\int_{\Omega\times \Theta}[(X(w)\otimes ( g'(\theta)+\ell'_w(\theta) )]\lambda(\theta)e^{\ell_w(\theta)}\nu(dw)d\theta=0.\end{equation}
From (\ref{N}) and Stokes we get that $$\int _{\Theta}[\theta\otimes(g'(\theta)+\ell'_w(\theta))]\lambda(\theta)e^{\ell_w(\theta)}d\theta=-I_k\int _{\Theta} \lambda(\theta)e^{\ell_w(\theta)}d\theta$$ and thus integrating with respect to $\nu(dw)$ we have 
\begin{equation}\label{N'}\int_{\Omega\times \Theta}[\theta\otimes ( g'(\theta)+\ell'_w(\theta) )]\lambda(\theta)e^{\ell_w(\theta)}\nu(dw)d\theta=-I_k.\end{equation}
since $P(d\theta,dw)=\lambda(\theta)e^{\ell_w(\theta)}\nu(dw)d\theta$ is a probability measure on $\Omega\times \Theta.$ 
We now combine (\ref{M'}) and (\ref{N'}) to get finally 
$$\int_{\Omega\times \Theta}[(X(w)-\theta)\otimes ( g'(\theta)+\ell'_w(\theta) )]\lambda(\theta)e^{\ell_w(\theta)}\nu(dw)d\theta=I_k.$$ Now consider the expectation under $ P(d\theta,dw)$ of the random variable 
$$(\theta,w)\mapsto [(X(w)-\theta),( g'(\theta)+\ell'_w(\theta) ]\otimes [(X(w)-\theta),( g'(\theta)+\ell'_w(\theta) ]$$
which is valued in the set of semi positive definite matrices of order $2k.$ This expectation is nothing but the matrix 
(\ref{vT}) which is therefore semi positive definite also. To prove this, the only point left is to check that 
$$\int_{\Omega\times \Theta}[ ( g'(\theta)+\ell_w(\theta) )  \otimes ( g'(\theta)+\ell_w(\theta) )]\lambda(\theta)e^{\ell_w(\theta)}\nu(dw)d\theta=I_{\lambda}+\int_{ \Theta}I(\theta)\lambda(\theta)d\theta.$$ This comes from 
$$\int_{\Omega}\ell'_w(\theta)e^{\ell_w(\theta)}\nu(dw)=0\Rightarrow \int_{\Omega}[g'(\theta)\otimes \ell'_w(\theta)]e^{\ell_w(\theta)}\nu(dw)=0\Rightarrow \int_{\Omega\times \Theta}[g'(\theta)\otimes \ell'_w(\theta)]P(d\theta,dw)=0.$$ \hfill
$ \square$

\vspace{4mm}\noindent \textbf{Remarks.} Comparison with the Bayesian Cram\'er-Rao inequality is in order. This inequality says that if $X$ is an unbiased estimator of $\theta$ then 

\begin{equation}\label{BCR}\int_{\Theta}(X-\theta)\otimes (X-\theta))\lambda(\theta)d\theta-\left(\int_{\Theta} I(\theta) \lambda(\theta)d\theta)\right)^{-1}\end{equation} is semi positive definite.

The  Van Trees inequality is about possibly \textit{biased estimators} while the inequality \eqref{BCR}  is about the unbiased case. Applying the Van Trees inequality to the unbiased case provides just a weaker result since saying by \eqref{BCR} that 

$$A=\left[\begin{array}{cc}\int_{\Omega\times \Theta}[(X(w)-\theta)\otimes (X(w)-\theta)]P_{\theta}(dw)\lambda(\theta)d\theta&I_k\\I_k&\int_{ \Theta}I(\theta)\lambda(\theta)d\theta\end{array}\right]$$ is semi positive definite is stronger than  saying that  
$A+\left[\begin{array}{cc}0&0\\0&I_{\lambda}\end{array}\right]$ is semi positive definite.

On the other hand, introduce  the expectation $\psi(\theta)$ of the biased estimator $X$ of $\theta$ . Consider the following semi positive definite matrix of order $k$ $$C=\int_{ \Theta}[(\psi(\theta)-\theta)\otimes (\psi(\theta)-\theta)]\lambda(\theta)d\theta.$$ 
Thus we have 
$$\int_{\Omega\times \Theta}[(X(w)-\theta)\otimes (X(w)-\theta)]P_{\theta}(dw)\lambda(\theta)d\theta=C+\int_{\Theta}\mathrm{Cov}_{\theta}(X)\lambda(\theta)d\theta.$$ The general Bayesian Cram\'er-Rao inequality  says that the matrix 
$$B=\left[\begin{array}{cc}\int_{ \Theta}\mathrm{Cov}_{\theta}(X)\lambda(\theta)d\theta& \int_{ \Theta}\psi'(\theta)\lambda(\theta)d\theta \\ \int_{ \Theta}\psi'(\theta)^t\lambda(\theta)d\theta&\int_{ \Theta}I(\theta)\lambda(\theta)d\theta\end{array}\right]$$ is semi positive definite 
and thus trivially if $B_1=\left[\begin{array}{cc}C&0\\0&I_{\lambda}\end{array}\right]$ then 
$B+B_1$ is semi positive definite. The  Van Trees inequality rather says  that if we consider  $$B_2=\left[\begin{array}{cc}C&I_k-\int_{ \Theta}\psi'(\theta)\lambda(\theta)d\theta\\I_k-\int_{ \Theta}\psi'(\theta)^t\lambda(\theta)d\theta&I_{\lambda}\end{array}\right]$$ then 
$B+B_2$ is semi positive definite. \textit{This fact is not trivial} since  $B_2$ is not positive definite in general. 
\section{Acknowledgements}I am indebted to Suhasini Subba Rao for patiently explaining to me the importance of dealing with biased estimators for the Van Trees inequality and thus the difference between $B_1$  and $B_2.$  My thanks  also go to Jean-Yves Tourneret who has introduced me to the general  problems of this note through the paper  \cite{BESSON2008} devoted to analogous questions for complex Gaussian laws, and to Christian Genest who has greatly improved the presentation of the results.

\end{document}